\documentclass[12pt]{article}

\usepackage{amssymb}
\usepackage{latexsym}
\def\part#1{\frac{\partial\phantom{q}}{\partial#1}}

\newenvironment{remark}{\begin{trivlist}\item[]{\bf Remark:} }
{\end{trivlist}}
\newenvironment{example}{\begin{trivlist}\item[]{\bf Example:} }
{\end{trivlist}}
\newenvironment{proof}{\begin{trivlist}\item[]{\bf Proof:} }
{\hfill $\Box$ \end{trivlist}}

\newtheorem{theorem}{Theorem}

\newtheorem{lemma}[theorem]{Lemma}

\def\tr{\mathop{\rm tr}\nolimits}

\def\Diff{\mathop{\rm Diff}\nolimits}

\begin{document}
\title{Stable forms and special metrics}
 \author{Nigel Hitchin\\[5pt]
\itshape  Mathematical Institute\\
\itshape 24-29 St Giles\\
\itshape Oxford OX1 3LB\\
\itshape UK\\
 hitchin@maths.ox.ac.uk}
\maketitle

\centerline{\textit{Dedicated to the memory of Alfred Gray.}}

\vskip .5cm

\begin{abstract}
We show how certain diffeomorphism-invariant functionals on differential forms in dimensions $6,7$ and $8$ generate in a natural way  special
geometrical structures in these dimensions: metrics of holonomy
$G_2$ and $Spin(7)$, metrics with weak holonomy $SU(3)$ and $G_2$,
and a new and unexplored example in dimension $8$. The general
formalism becomes a practical tool for calculating homogeneous or
cohomogeneity one examples, and we illustrate this with some newly discovered examples of $Spin(7)$ and $G_2$ metrics.
\end{abstract}

\section{Introduction}
 One of Alfred Gray's  most original concepts was that of
weak holonomy \cite{G}. This was an idea  clearly ahead of its time, as became evident in the
later work on Killing spinors of Baum et al. \cite{B}, setting it in a  natural context.  We shall
give here another natural approach to both weak holonomy and
special holonomy in low dimensions through the use of certain
invariant functionals of differential forms. This provides both a
canonical setting for these structures and sometimes an effective
means of finding them.

Our starting point is the question: ``What is a non-degenerate
form?". Symplectic geometry is the geometry determined by a closed non-degenerate $2$-form where  non-degeneracy  means
$\omega^m\ne0$. Another way to describe such
genericity at each point is to note that the orbit of a nondegenerate $\omega\in\Lambda^2V^*$ under the natural
action of $GL(V)$ on $\Lambda^2V^*$ is open. We shall say in general that
$\rho\in\Lambda^pV^*$ is  \textit{stable} if it lies  in an open
orbit. (The use of this word rather than ``non-degenerate" avoids  possible confusion: a $2$-form in odd dimensions may be stable but is always degenerate as a
bilinear form. Also, $\rho$ is stable in the sense of deformation
invariance: all forms in a neighbourhood of $\rho$ are
$GL(V)$-equivalent to $\rho$).

When does stability occur? Clearly not very often since $\dim
GL(V)=n^2$ is usually much smaller than $\dim \Lambda^p
V^*=n!/p!(n-p)!$. The result (over the complex numbers) is
classical \cite{R}, \cite{S}, \cite{Gu1}, \cite{Gu2} but Robert Bryant explained to me
all the real cases. Apart from the obvious case of $p=1,2$ there
are essentially only three more, where $p=3$ and $n=6,7,8$. In
these three cases the stabilizer subgroup of $\rho \in
\Lambda^p V^*$ in $GL(V)$  is a real form of one of the complex groups $$SL(3)\times
SL(3),\quad G_2,\quad PSL(3)$$ respectively. We shall be concerned 
here only with the real forms $SL(3,\mathbf{C})$ and the compact
groups $G_2$ and $PSU(3).$ Thus if a manifold $M^n$ admits a global stable $p$-form, it has a $G$-structure where $G$ is one of these groups.

Note that if $GL(V)$ has an open orbit
in $\Lambda^p V^*$ then it also does on the dual space $\Lambda^p
V\cong \Lambda^{n-p} V^*\otimes\Lambda^n V$. Since for $p\ne 0$
the scalars act non-trivially, there is then an open orbit on
$\Lambda^{n-p} V^*$. Thus if we are in a dimension $n$ where stable $p$-forms exist,
we can also consider stable $(n-p)$-forms. \vskip .25cm

The three stabilizers of $\rho$ above, as well as the symplectic
group, each preserve a volume element $\phi(\rho)\in \Lambda^nV^*$. Thus if we
have a compact oriented manifold $M^n$ and a $p$-form $\rho$ which
is everywhere stable, we can integrate $\phi(\rho)$ to obtain a
volume $V(\rho)$. Openness of the orbit implies that nearby forms
are also stable, so that the volume functional is defined and
smooth on an open set of forms.  \vskip .25cm

We now set up a number of variational problems involving this functional. First  we consider  a  critical point of $V(\rho)$ restricted to a
fixed cohomology class of closed forms in $H^p(M,{\bf R})$, performing a non-linear version of Hodge theory. If
$\rho$ is everywhere stable we find the following structures:

\begin{itemize}
\item
 for
$n=2m$ and $p=2$ or $2m-2$, a symplectic manifold, 

\item for $n=6$
the structure of a complex 3-manifold with trivial canonical
bundle (this is described in some detail in \cite{H}),

 \item for $n=7$ and $p=3$ or
$4$ a Riemannian manifold with holonomy $G_2$.
\end{itemize}
\vskip .25cm The critical points for $n=8$ and $p=3$ or $5$ form a
class of geometric structures which is  largely unexplored (though
I have benefited from Robert Bryant's thoughts on these). We study
them briefly in Section \ref{eight}, showing that they admit a
solution to the Rarita-Schwinger equations, which imposes
constraints on the Ricci tensor.  The lack of concrete compact examples
beyond the $8$-manifold $SU(3)$ itself is currently a stumbling
block in taking the analysis of these further, but the fact that
they arise from the same variational origins as $G_2$-manifolds
suggests that they ought to exist in abundance. \vskip .25cm
  To continue with other variational characterizations, we note that there is a canonically
defined indefinite quadratic form $$Q(d\alpha)=\int_M\alpha\wedge
d\alpha$$ on the space of exact $4$-forms $d\alpha$ on a
$7$-manifold $M$. Thus $Q$ defines an indefinite metric
on the closed $4$-forms in a fixed cohomology class. We then find:

\begin{itemize} \item
 the stable critical points of
$V(\rho)$ on the trivial cohomology class, subject to the
constraint $Q(d\alpha)=const$ define a $7$-manifold with a weak
holonomy $G_2$ structure, 

\item
   the
gradient flow of $V(\rho)$ on a fixed degree $4$ cohomology class
 yields a Riemannian metric of holonomy $Spin(7)$ on $M\times
{\mathbf R}$.
\end{itemize}
\vskip .25cm
 Finally, consider a $6$-manifold $M$ and a stable closed
$3$-form $\rho$ together with a stable closed $4$-form $\sigma$. We note
that the spaces of exact $3$-forms $d\alpha$ and exact $4$-forms
$d\beta$  on a $6$-manifold are formally dual to each other via
the pairing $$\langle d\alpha,d\beta\rangle=\int_M\alpha\wedge
d\beta.$$ We find

\begin{itemize}
\item
for $\rho$ and $\sigma$ exact, the stable critical points of
$3V(\rho)+8V(\sigma)$ subject to the constraint $\langle
\sigma,\rho\rangle=const$ define a manifold with weak holonomy
$SU(3)$,

\item the pairing defines a  formal symplectic structure on the
product of a degree $3$ and degree $4$ cohomology class and the
Hamiltonian flow of the functional $V(\rho)-2V(\sigma)$, with an
appropriate initial condition, generates a Riemannian metric with
holonomy $G_2$ on $M\times {\mathbf R}$.
\end{itemize}
(The particular coefficients here are a computational convenience and have no geometrical meaning.)
\vskip .25cm
 The formal variational setting of all these special
geometries is perfectly general -- any such structure appears in
this way -- but it  also has a practical value. In the homogeneous
or cohomogeneity one situation we merely set up the same
variational equations using \textit{invariant} forms and easily
derive the corresponding equations. We demonstrate this in the
case of $S^7$ by deriving the equations recently used
by \v Cvetic et al. \cite{C} to find a new $Spin(7)$ metric on
$\mathbf{R}^8$, and in the work of Brandhuber et al. \cite{Br} for a new example of a metric of holonomy $G_2$.

The author wishes to thank the Universidad Aut\'onoma, Madrid
 and the Programa Cat\`edra Fundaci\'on Banco de Bilbao y Vizcaya  for support during the preparation of this
paper.

\section{The linear algebra of stable forms}

Let $\rho \in \Lambda^p V^*$ be {\textit stable}, in the sense
described above, i.e. it lies in an open orbit $U$ of $GL(V)$. We
consider the cases:

\begin{enumerate}
\item
$\dim V=2m$, $p=2$ or $2m-2$: stabilizer $Sp(2m,\mathbf{R})$
\item
$\dim V=6$, $p=3$: stabilizer $SL(3,\mathbf{C})$
\item
$\dim V=7$, $p=3$ or $4$: stabilizer $G_2$
\item
$\dim V=8$, $p=3$ or $5$: stabilizer $PSU(3)$
\end{enumerate}

We see that each stabilizer preserves a volume form: $G_2$ and
$PSU(3)$ are compact so this is the volume form of an invariant positive
definite inner product on $V$. The symplectic group fixes the
Liouville volume and $SL(3,\mathbf{C})$ preserves a complex
$3$-form $\Omega$ and hence the real $6$-form
$i\Omega\wedge\bar\Omega$. In the appendix we
give the concrete expression of this form in each case. It is algebraically
determined by $\rho$ and smooth on $U$.
\vskip.25cm

 The volume form $\phi(\rho)$ associated to $\rho$ defines a
$GL(V)$-invariant map $$\phi:U\rightarrow \Lambda^ nV^*.$$ Applying invariance to the action of the scalar matrices, we see that  $$\phi(\lambda^p\rho)=\lambda^n\phi(\rho)$$ so
that $\phi$ is homogeneous of degree $n/p$.
\vskip .25cm
 The
derivative of $\phi$ at $\rho$ is an invariantly defined element
of $(\Lambda^pV^*)^*\otimes \Lambda^nV^*$. Since
$(\Lambda^pV^*)^*\cong \Lambda^{n-p}V^*\otimes \Lambda^nV$, the derivative
lies in $\Lambda^{n-p}V^*$ so there is a unique element
$\hat\rho\in \Lambda^{n-p}V^*$ for which
\begin{equation}\label{deriv}
D\phi(\dot\rho)=\hat\rho\wedge\dot\rho \end{equation} Taking
$\dot\rho=\rho$, Euler's formula for a homogeneous function gives
\begin{equation}\label{hom}
\hat\rho\wedge\rho=\frac{n}{p}\phi(\rho)
\end{equation}
\vskip .25cm
\begin{example}
If $\dim V=2m$ and $\omega\in \Lambda^2V^*$ is stable (meaning
non-degenerate here), then we take the Liouville volume form $$\phi(\omega)=\frac{1}{m!}\omega^m\in
\Lambda^{2m}V^*$$ This is clearly homogeneous in $\omega$ of
degree $m=2m/2=n/p$. Differentiating, we see that
$$D\phi(\dot\omega)=\frac{1}{(m-1)!}\dot\omega\wedge\omega^{m-1}$$
so that $$\hat\omega=\frac{1}{(m-1)!}\omega^{m-1}.$$

We may also consider a stable  $\rho \in \Lambda^{2m-2}V^*$, which
is in the open orbit $U$ consisting of forms
$\rho=\omega^m/(m-1)!$ for a non-degenerate $\omega$. In this case
we have $\hat\rho=\omega/(m-1)\in \Lambda^2V^*$.
\end{example}
\vskip .25cm
The precise form of $\hat\rho$ is determined by seeing
which elements of $\Lambda^{n-p}V^*$ are fixed by the stabilizer.
The symplectic case is done in the example above. For the others
we see easily that:

\begin{itemize}
\item
for $n=6$, $p=3$, $\hat\rho$ is determined by the property that
$\Omega=\rho+i\hat\rho$ is a complex $(3,0)$-form preserved by
$SL(3,\mathbf{C})$,
\item
for $n=7$, $p=3$ or $4$, $\hat\rho=\ast\rho$, where $\ast$ is the
Hodge star operator for the inner product on $V$
\item
for $n=8$, $p=3$ or $p=5$, $\hat\rho=- \ast\rho$.
\end{itemize}

\begin{remark} There is clearly a choice in what we call the volume in each case. There are conventions -- the Liouville volume $\omega^m/m!$ in the symplectic case for example. For most purposes it makes no difference to the results that follow, but when we need to find a metric, as in Sections $5$ and $6$, we shall make a more explicit choice to aid the calculations.
\end{remark}

\section{Critical points}

Suppose now that $M$ is a closed, oriented $n$-manifold. If
$\rho\in \Omega^p(M)$ is a global $p$-form, then it is a section
of $\Lambda^pT^*$. Suppose its value is stable at each point (such an assumption requires of course the reduction of the structure group of the tangent bundle to one of the stabilizers above). Then
we can define a  functional by taking the total volume:
$$V(\rho)=\int_M\phi(\rho).$$ By definition of stability, nearby
forms will be stable and so we can differentiate the functional. We
shall set up a variational problem by considering the volume
restricted to \textit{closed} stable $p$-forms in a given
cohomology class.

\begin{theorem}\label{var}
A closed stable form $\rho\in \Omega^p(M)$ is a critical point of
$V(\rho)$ in its cohomology class if and only if $d\hat\rho=0.$
\end{theorem}
\begin{proof}
Take the first variation of $V(\rho)$: $$\delta V(\dot\rho)=\int_M
D\phi(\dot\rho)=\int_M\hat\rho\wedge\dot\rho$$ from (\ref{deriv}).
But the variation is within a fixed cohomology class so
$\dot\rho=d\alpha$. Thus $$\delta V(\dot\rho)=\int_M\hat\rho\wedge
d\alpha=\pm\int_Md\hat\rho\wedge \alpha$$ and the  variation
vanishes for all $d\alpha$ if and only if $$d\hat\rho=0.$$
\end{proof}
\begin{example}
A closed stable $2$-form $\omega$ on an even-dimensional manifold is a
symplectic form. The volume $V(\omega)$ is then \textit{constant}
on a fixed cohomology class $[\omega]$ since it is just the
evaluation of the cup product:
$$V(\omega)=\frac{1}{m!}[\omega]^m[M].$$ Rather trivially a
symplectic manifold appears as a critical point here.

On the other hand if we take $\rho=\omega^{m-1}/(m-1)!$ to be the closed
form then the condition $d\rho=0$ is (for $m>2$) weaker than
$d\omega=0$. The functional genuinely varies and the critical
points are where $$(m-1)d\hat\rho=d\omega=0.$$ This is an alternative
way of obtaining a symplectic manifold as a critical point.
\end{example}

 \begin{example} A $7$-manifold $M$ which has either a closed stable $3$-form or a
closed stable $4$-form which is a critical point for $V$ in its
cohomology class has the structure of a Riemannian manifold with
holonomy $G_2$. This follows from the theorem of M.~Fern\'andez
and A.~Gray \cite{F} that the holonomy reduces to $G_2$  if and
only if $$d\rho=d\ast\rho=0.$$ \end{example}

 \begin{example} The case of $n=6$ means that we have a complex closed locally decomposable $3$-form $\rho+i\hat\rho$. This
yields the structure of a complex $3$-manifold with trivial canonical bundle
 as shown in \cite{H}.
 \end{example}
 The $8$-dimensional case will occupy us
 next.

 \section{Eight-manifolds with $PSU(3)$ structure}\label{eight}

Suppose that $M$ is a compact $8$-manifold with a stable $3$-form
$\rho$ such that, with respect to the metric determined by $\rho$,  $$d\rho=d\ast\rho=0.$$ An example is $SU(3)$
itself, where $\rho$ is covariant constant and is a multiple of the
standard bi-invariant form $$\rho=\tr(g^{-1}dg)^3.$$

We shall never find a metric on a compact simply-connected $M^8$ whose
holonomy is $PSU(3)$ other than this example, from Berger's
classification of Riemannian holonomy groups. This is not the case of weak holonomy either --
there are no Killing spinors. We do however have one interesting
object: a Rarita-Schwinger field. In  physics terminology this
is a spin $3/2$ field. We take one of the two  spinor bundles (say $S^{+}$)
and consider a spinor-valued $1$-form -- a section $\gamma$ of $S^+\otimes \Lambda^1$. This satisfies
the \textit{Rarita-Schwinger equation} if $$D\gamma=0$$ and
$$d^*\gamma=0.$$ Here $$D:C^{\infty}(S^+\otimes
\Lambda^1)\rightarrow C^{\infty}(S^-\otimes \Lambda^1)$$ is the
Dirac operator with coefficients in the bundle of $1$-forms
$\Lambda^1$ and $$d^*:C^{\infty}(S^+\otimes \Lambda^1)\rightarrow C^{\infty}(S^+)$$ is the covariant $d^*$ operator on $1$-forms with
coefficients in the spinor bundle $S^+$.

\vskip .25cm

First we
describe $\gamma$ as a $PSU(3)$-invariant object.

\begin{lemma} Let $S^+,S^-$  be the two spin
representations and $\Lambda^1$ the standard vector representation
of $Spin(8)$. Then restricted to the lift $PSU(3)\subset Spin(8)$
of the adjoint representation, these three
representations are equivalent.
\end{lemma}
\begin{proof} If $\pm x_1,\dots,\pm x_4$ are the weights of the
$8$-dimensional vector representation of $Spin(8)$, the weights
of the spin representations $S^{\pm}$ are
\begin{equation}
\frac{1}{2}(\pm x_1\pm
x_2\dots\pm x_4)
\label{spinweights}
\end{equation} where there is an even number of minus signs for
$S^+$ and an odd number for $S^-$.

If $\alpha,\beta,\alpha+\beta$ are the positive roots of $SU(3)$
then substituting $$x_1=0,\quad x_2=\alpha,\quad x_3=\beta,\quad
x_4=\alpha+\beta$$ we have from (\ref{spinweights}) for $S^+$ and $S^-$  the same weights $0,\pm \alpha,\pm\beta,\pm \alpha+\beta$ as the adjoint
representation.
\end{proof}

\begin{remark} The lemma implies a rather interesting property of
the Lie algebra of $SU(3)$. Recall that Clifford multiplication of
vectors on spinors is skew adjoint and satisfies $x^2=-(x,x)1$. In
eight dimensions the three representations $S^+$, $S^-$ and
$\Lambda^1$ are all real and $8$-dimensional (this is triality).
If $\varphi\in S^{\pm}$ and $x\in \Lambda^1$ then
$$(x\varphi,x\varphi)=(-x^2\varphi,\varphi)=(x,x)(\varphi,\varphi)$$
so that  Clifford multiplication $$\Lambda^1\otimes S^+\rightarrow
S^-$$ is an orthogonal product. From the lemma all three
representations are equivalent under $SU(3)$ so we must have an
$SU(3)$-invariant  orthogonal multiplication on the Lie algebra of
$3\times 3$ skew-hermitian matrices with trace zero. Here it is:
$$A\times B=\omega AB-\bar\omega
BA-\frac{i}{\sqrt{3}}\tr(AB)I$$ where
$\omega=(1+i\sqrt{3})/2$.
\end{remark}
\vskip .25cm The $SU(3)$-invariant isomorphism $S^+\cong \Lambda^1$
defines a section of $(S^+)^*\otimes \Lambda^1\cong S^+\otimes
\Lambda^1$ and this is  what we take to be $\gamma$. There is of course an equivalent
section of $S^-\otimes \Lambda^1$. We now prove:

\begin{theorem} If $M$ is an $8$-manifold with a $PSU(3)$ structure defined by
a $3$-form $\rho$ with $d\rho=d\ast\rho=0$, then the section
$\gamma \in C^{\infty}(S^+\otimes \Lambda^1)$ satisfies the
Rarita-Schwinger equations.
\end{theorem}
\begin{proof} The covariant derivative of the $3$-form $\rho$ at any point can be written
$$\nabla \rho=A(\rho)$$
where $A\in \Lambda^2\otimes \Lambda^1$ and acts on $\rho$ by identifying $\Lambda^2$ with the Lie algebra of $SO(8)$. Since $\rho$ is fixed by $SU(3)$ we may as well assume that $A\in \Lambda_0^2\otimes \Lambda^1$ where $\Lambda^2_0$ is the orthogonal complement in $\Lambda^2$ to the Lie algebra of $SU(3)\cong\Lambda^1$. Since $\rho$ is not covariant constant, $A$ will not vanish. However, the harmonicity condition on $\rho$ will force many of its components (as representations of $SU(3)$) to vanish.

We shall index representations by their highest weight -- the adjoint representation $\Lambda^1$ of $SU(3)$ has highest weight $\alpha +\beta$. If we   decompose the tensor product $\Lambda^1\otimes\Lambda^1$ into irreducible representations as in Chapter 6 of Salamon's book \cite{Sal}, we find that $\Lambda_0^2$ is the direct sum of two $10$-dimensional irreducibles $\Lambda^2_+,\Lambda^2_-$ with highest weight $2\alpha+\beta$ and $\alpha+2\beta$ respectively. These are interchanged under a change of orientation, just like self-dual and anti-self-dual forms in four dimensions. Similarly $V(2\alpha+\beta)\otimes V(\alpha +\beta)$ breaks up into  irreducibles with highest weights
$$3\alpha+2\beta,\quad 2\alpha+2\beta,\quad 2\alpha+\beta,\quad \alpha+\beta$$
and each with multiplicity one.
Robert Bryant informed me of the following result:

\begin{lemma} If $d\rho=d\ast\rho=0$, the components of $A$ with highest weight $\alpha+\beta$, $2\alpha+\beta$, $\alpha+2\beta$ and $2\alpha+2\beta$ all vanish.
\end{lemma}
\begin{proof}
Since we are considering exterior powers of the Lie algebra it is convenient to think of these as spaces of left-invariant forms on the group. Then, for example, $\Lambda^2_0$ is the space of coclosed invariant $2$-forms. The $3$-form $\rho$ is built out of the structure constants of the Lie algebra  $\mathfrak {su}(3)$ and as a consequence of this, the action of $a\in\Lambda^2$ on $\rho$ can be rewritten as $da$ if we consider $a$ as a $2$-form on the group.  Thus, if
$$A=\sum_i a_i\otimes e_i\in \Lambda_0^2\otimes \Lambda^1$$
then
$$d\rho=\sum_i e_i\wedge da_i$$ and
$$d^*\rho=\sum_i \iota(e_i)da_i.$$
Now if $x_{\gamma}\in \mathfrak {su}(3)\otimes \mathbf{C}$ denotes a root vector for the root $\gamma$, the vector
$$A=x_{\alpha}\wedge x_{\alpha+\beta}\otimes x_{\beta}$$ lies in $V(2\alpha+2\beta)\subset V(2\alpha+\beta)\otimes V(\alpha +\beta)$. Let $h$ be an element of the Cartan subalgebra, then because $\iota(h)dx_{\alpha}=\mathcal{L}_hx_{\alpha}=\alpha(h)x_{\alpha}$ etc.
$$\iota(h)d\rho=\iota(h)(x_{\beta}\wedge d(x_{\alpha}\wedge x_{\alpha+\beta}))=\beta(h)x_{\beta}\wedge x_{\alpha+\beta}\wedge x_{\alpha}.$$
Since we can find $h$ for which $\beta(h)\ne 0$ this shows that the irreducible representation $V(2\alpha+2\beta)\subset \Lambda^2_+\otimes \Lambda^1$ maps non-trivially into $\Lambda^3$ under $A(\rho)$. That representation is also contained in $\Lambda^2_-\otimes \Lambda^1$, interchanging the roles of $\alpha$ and $\beta$. However, since we can choose $h$ such that $\alpha(h)=0$ and $\beta(h)\ne 0$ and vice-versa, $V(2\alpha+2\beta)$ appears twice in $\Lambda^3$, and so if $d\rho=0$, then both of these components in $A$ vanish.
\vskip .25cm
We now work similarly with $V(\alpha+\beta)\subset V(2\alpha+\beta)\otimes V(\alpha+\beta)$ using the vector
$$A= x_{\alpha}\wedge x_{\alpha+\beta}\otimes x_{-\alpha}.$$
We deduce that if  $d\rho=0$, then both of the $V(\alpha+\beta)$ components in $A$ vanish.
\vskip .25cm
Now consider $V(2\alpha+\beta)=\Lambda^2_+\subset \Lambda^2_+\otimes \Lambda^1$. The inclusion is unique and is given by
$$a\mapsto\sum_i\mathcal {L}_{e_i}a\otimes e_i$$
for an orthonormal basis $\{e_i\}$ of the Lie algebra $\mathfrak{su}(3)\cong\Lambda^1$. In this case
$$d^*\rho=\sum_i\iota(e_i)d(\mathcal{L}_{e_i}a)=\sum_i\mathcal{L}^2_{e_i}a-d(\iota(e_i)\mathcal{L}_{e_i} a)$$
and $a\mapsto \iota(e_i)\mathcal{L}_{e_i}a$ is an invariant map from the irreducible  $\Lambda^2_+$ to the irreducible $\Lambda^1$ and thus must vanish. Hence
$d^*\rho$ is the Casimir $\sum_i\mathcal{L}^2_{e_i}a$ which is a non-zero scalar multiple of $a$. It follows that if $d^*\rho=0$, the component $V(2\alpha +\beta)$ (and similarly $V(\alpha+2\beta)$) in $A$ vanishes. This completes the proof of the lemma.
\end{proof}
To return to the theorem, consider $D\gamma\in C^{\infty}(S^-\otimes \Lambda^1)$. The right hand side, since $S^-\cong\Lambda^1$, is a vector bundle associated to the representation $V(\alpha +\beta)\otimes V(\alpha+\beta)$. The skew part we calculated to have highest weights $2\alpha+\beta,\alpha+2\beta$ and $\alpha+\beta$. The calculation gives at the same time the symmetric part to have weights $0$, $\alpha+\beta$ and $2\alpha+2\beta$. But $D\gamma$ is the image under an invariant map of $A$ which from the lemma only has components with highest weights $3\alpha + 2\beta$ and $2\alpha+3\beta$. Since these do not occur in $S^-\otimes \Lambda^1$ we deduce that $D\gamma=0$. Similarly, since $S^+=V(\alpha +\beta)$, we see that $d^*\gamma=0$.
\end{proof}
 The existence of a Rarita-Schwinger field implies the vanishing of certain components of the Ricci tensor. In a sense the Einstein equations are the integrability condition for a linear system given by the Rarita-Schwinger operator \cite{J}. More concretely, given $\gamma\in C^{\infty}(S^+\otimes \Lambda^1)$ we form the second covariant derivative
$$\nabla^2\gamma\in C^{\infty}(S^+\otimes\Lambda^1\otimes  \Lambda^1\otimes\Lambda^1)$$
with components $\gamma_{i;jk}$ in an orthonormal basis. Covariantly differentiating the  equation $D\gamma=0$ gives
$$\sum_j e_j\gamma_{i;jk}=0$$
and by contraction
\begin{equation}
\sum_{i,j} e_j\gamma_{i;ji}=0.
\label{ji}
\end{equation}
But differentiating the equation $d^*\gamma=0$ gives
$$\sum_{i,j}\gamma_{i;ij}=0$$
and hence also
\begin{equation}
\sum_{i,j} e_j\gamma_{i;ij}=0.
\label{ij}
\end{equation}
The difference of (\ref{ji}) and (\ref{ij}) is skew-symmetric in the last two indices and can therefore be rewritten in terms of the action of the curvature tensor on $S^+\otimes \Lambda^1$. This only involves the Ricci tensor:
\begin{equation}
\sum_{i,j}R_{ij}e_i\gamma_j=0
\label{ricci}
\end{equation}
The expression $\sum_{i,j}R_{ij}e_i\gamma_j$ in (\ref{ricci}) defines an $SU(3)$-invariant map from the space of
symmetric tensors $R_{ij}$ to $S^-$. We have seen that $S^-\cong V(\alpha+\beta)$ and this representation occurs with multiplicity one in $Sym^2 \Lambda^1$. Thus equation (\ref{ricci}) implies the vanishing of $8$ of the $36$ components of the Ricci tensor.

\begin{remark} One consequence of the theorem is that the operator
$$\mathcal {D}:\Omega^1(M)\rightarrow\Omega^1(M)$$
defined by the orthogonal multiplication 
$$\mathcal {D}\alpha=\sum_i e_i\times\nabla_i\alpha$$
is equivalent (under the isomorphisms $\gamma$ for $S^+$ and $S^-$) to the Dirac operator. The  multiplication $\alpha\times\beta$ is of course Clifford multiplication under these isomorphisms, but the Levi-Civita connections on $\Lambda^1$  ond on $S^{\pm}$ are different. Nevertheless,  $D\gamma=0$ implies that the Dirac operators correspond.
\end{remark}

\section{Constrained critical points}

On a compact oriented manifold $M^n$ there is a non-degenerate pairing between the spaces of forms $\Omega^p(M)$ and $\Omega^{n-p}(M)$ defined by
$$\int_M\alpha\wedge\beta.$$
If $\alpha=d\gamma\in \Omega^p(M)$ is exact,  then by Stokes'theorem
$$\int_M d\gamma\wedge\beta=(-1)^{p}\int_M\gamma\wedge d\beta$$
which vanishes for all $\gamma$ if and only if $\beta$ is closed. We thus have a non-degenerate pairing between $\Omega^p_{exact}(M)$ and
$$\Omega^{n-p}(M)/\Omega^{n-p}_{closed}(M)$$
Since the exterior derivative $d$ maps this latter space isomorphically onto $\Omega^{n-p+1}_{exact}(M)$, then formally (in the sense of non-degenerate pairings) we can say that
\begin{equation}
\Omega^p_{exact}(M)^*\cong \Omega_{exact}^{n-p+1}(M)
\label{duality}
\end{equation}

Now consider $n=7$ and $p=4$. We have seen that a metric of holonomy $G_2$ arises from a stable critical point for the functional $V(\rho)$ on a fixed cohomology class  $[\rho]$ of closed $4$-forms. On the other hand since $\ast\rho\wedge \rho=7\phi(\rho)/4$ and $\rho$ is both closed and coclosed, we must have
$[\ast\rho]\cup[\rho][M]=V(\rho)\ne 0$ and in particular $[\rho]\ne 0$.
For the trivial cohomology class it follows that there are  no stable critical points for this functional. However, we have from (\ref{duality})
$$\Omega^4_{exact}(M)^*\cong \Omega_{exact}^{4}(M)$$
and hence a non-degenerate quadratic form on $\Omega_{exact}^{4}(M)$ given by
$$Q(d\gamma)=\int_M\gamma\wedge d\gamma.$$
We now have two natural functionals on the space of exact $4$-forms: $V$ and $Q$.
\begin{theorem}\label{weakg2}
 An  exact stable $4$-form $\rho$ on a compact $7$-manifold is a critical point of $V(\rho)$ subject to the constraint $Q(\rho)=const.$ if and only if $\rho$ defines a metric with weak holonomy $G_2$.
\end{theorem}
\begin{proof}
From Theorem (\ref{var}), the first variation of $V$ at $\rho=d\gamma$ is
$$\delta V(d\dot\gamma)=\int_M\ast \rho\wedge d\dot\gamma$$
and the first variation of the quadratic form $Q$ is
$$\delta Q(d\dot\gamma)=2\int_M \dot\gamma\wedge\rho.$$
Thus, introducing a Lagrangian multiplier, the constrained critical point is given by
$$d(\ast \rho)=\lambda\rho$$
and from \cite{Fr}, this is equivalent to the structure of a manifold with weak holonomy $G_2$ (sometimes called a \textit{nearly parallel} $G_2$ structure).
\end{proof}

\begin{remark} A $7$-manifold with weak holonomy $G_2$ is an Einstein manifold with positive scalar curvature. There are many examples, such as $3$-Sasakian manifolds (see \cite{Bo}) and their squashed versions (see  \cite{Fr}).
\end{remark}
\vskip .25cm
Next consider the case of $n=6$ and $p=3$. Then (\ref{duality}) tells us
$$\Omega^3_{exact}(M)^*\cong \Omega_{exact}^{4}(M)$$
so that the spaces of exact $3$-forms and exact $4$-forms are formally dual to each other. The pairing for an exact $3$-form $\rho=d\alpha$ and an exact $4$-form $\sigma=d\beta$ is
\begin{equation}
\langle \rho,\sigma\rangle=\int_M \alpha\wedge \sigma=-\int_M\rho\wedge \beta
\label{rspairing}
\end{equation}
We shall consider in $6$ dimensions a variational problem involving the three functionals $V(\rho),V(\sigma)$ and $\langle \rho,\sigma\rangle$, but before proving a theorem, let us see how $\rho$ and $\sigma$ define a reduction of the structure group of the tangent bundle of $M^6$ to $SU(3)$. The $3$-form $\rho$, being stable, provides a  reduction to $SL(3,\mathbf{C})$ and $\sigma=\omega^2/2$  to $Sp(6,\mathbf{R})$. The group $SU(3)$ is an intersection of these two but there are two compatibility conditions for $\rho$ and $\sigma$ to achieve this. The first is
\begin{equation}
\omega\wedge\rho=0
\label{compat1}
\end{equation}
From the point of view of $\rho$ and the complex structure it describes, this says that $\omega$ is of type $(1,1)$. Since $\rho+i\hat\rho$ is a $(3,0)$ form the property $\omega\wedge\hat\rho=0$ follows from (\ref{compat1}). From the point of view of the symplectic form $\omega$, the equation says that $\rho$ is primitive. For a stable $3$-form $\rho$ and a stable $4$-form $\sigma=\omega^2/2$, we shall say that the pair $(\rho,\sigma)$ is of positive type if the almost complex structure $I$ determined by $\rho$ (see  Appendix (8.2)) makes $\omega(X,IX)$ a positive definite form. This is clearly an open condition. If $\omega\wedge\rho=0$, the condition of positive type means that the Hermitian form defined by $\omega$ is positive definite.

The second condition is
\begin{equation}
\phi(\rho)=c\phi(\sigma)
\label{compat2}
\end{equation}
For a constant $c$. This says that the complex $3$-form $\rho+i\hat\rho$ has constant length  relative to the Hermitian metric.

Since $\rho$ and $\sigma$ satisfying these conditions define a metric, it is natural to normalize the constant $c$ above so that the volume forms $V(\rho)$ and $V(\sigma)$ are fixed multiples of the metric volume form. The most convenient way to do this is to see the normal form of the $G_2$ $3$-form $\varphi$ and its dual $\ast\varphi$ in terms of an orthonormal basis as in \cite{Sal}:
$$\varphi=e_7(e_5e_6+e_1e_4+e_3e_2)+(e_1e_2e_3-e_3e_4e_5+e_1e_3e_6-e_4e_2e_6)$$
$$\ast\varphi=e_7(e_3e_4e_6-e_1e_2e_6+e_1e_3e_5-e_4e_2e_5)+(e_1e_3e_4e_2+e_5e_6e_2e_3+e_5e_6e_4e_1).$$
The $6$-dimensional geometry is defined by  $\omega=e_5e_6+e_1e_4+e_3e_2$ and $\rho=e_1e_2e_3-e_3e_4e_5+e_1e_3e_6-e_4e_2e_6$, and  then $\hat\rho=e_3e_4e_6-e_1e_2e_6+e_1e_3e_5-e_4e_2e_5$. We  have from the above expressions
\begin{equation}
\phi(\sigma)=\frac{1}{6}\omega^3=\frac{1}{4}\hat\rho\wedge\rho=\frac{1}{2}\phi(\rho)
\label{volumes}
\end{equation}
 and
\begin{equation}
\varphi=e_7\wedge\omega+\rho,\quad\ast\varphi=e_7\wedge\hat\rho-\sigma
\label{ast}
\end{equation}
 We now prove the theorem.
\begin{theorem}\label{weaksu3}
 A pair $(\rho,\sigma)$ of exact, stable forms of positive type on a compact $6$-manifold forms a critical point of $3V(\rho)+8V(\sigma)$  subject to the constraint  $\langle\rho,\sigma\rangle=const.$ if and only if $\rho$ and $\sigma$ define a metric with weak holonomy $SU(3)$.
\end{theorem}
\begin{proof} From Theorem (\ref{var}), the first variation of $3V(\rho)+8V(\sigma)$ is
$$3\int_M\hat \rho\wedge \dot\rho+4\int_M\omega\wedge \dot\sigma.$$
 The first variation of $P=\langle\rho,\sigma\rangle$ is
$$\delta P(\dot\rho,\dot\sigma)=\int_M\dot\rho\wedge\beta+\int_M \dot\sigma\wedge\alpha$$
where $\rho=d\alpha$ and $\sigma=d\beta$.
Using Stokes' theorem and a Lagrange multiplier $12\lambda$ we find the constrained critical point to be given by the equations $3d\hat\rho=-12\lambda\sigma=-6\lambda\omega^2$ and $4d\omega=12\lambda\rho$, i.e. 
\begin{equation}
 d\hat\rho=-2\lambda\omega^2
\label{var1}
\end{equation}
\begin{equation}
d\omega=3\lambda\rho
\label{var2}
\end{equation}
The compatibility conditions (\ref{compat1}), (\ref{compat2}) actually follow from these equations. From (\ref{var2}) we have
$$\lambda\omega\wedge\rho= \frac{1}{3}\omega\wedge d\omega=\frac{1}{6}d(\omega^2)=0$$
since $\sigma=\omega^2/2$ is closed. Moreover, from (\ref{var1})
$$2\lambda\omega^3=-\omega \wedge d\hat\rho=-d(\omega\wedge\hat\rho)+d\omega\wedge\hat\rho=3\lambda\rho\wedge\hat\rho$$
using (\ref{var2}) and the fact that  $\omega\wedge\hat\rho=0$, which, as we have seen, follows from $\omega\wedge\rho=0$. Thus
$\omega^3=3\hat\rho\wedge\rho/2$ as in (\ref{volumes}).

The equations (\ref{var1}),(\ref{var2}) give a metric of weak holonomy $SU(3)$,  sometimes called a {\it nearly K\"ahler} metric. To see this, consider the $3$-form on $M\times \mathbf{R}$
$$\varphi=\frac{r^2}{\lambda}dr\wedge \omega+r^3\rho.$$
This is stable, and from (\ref{var2}) closed. Moreover, by comparing with the normal forms above, we see that the $G_2$ metric it defines is
$$(dr/\lambda)^2+r^2g$$
where $g$ is the $SU(3)$ metric. From (\ref{ast}) we have
$$\ast\varphi=\frac{r^3}{\lambda}dr\wedge\hat\rho-\frac{r^4}{2}\omega^2$$
and then from (\ref{var1}), $d\ast\varphi=0$. It follows that the cone metric above is a Riemannian metric of holonomy $G_2$. However, from \cite{Ba}, this implies that $M^6$ is nearly K\"ahler.
\end{proof}
\begin{remark} The particular coefficients of $V(\rho)$ and $V(\sigma)$ in the theorem are not crucial. Any two positive numbers will, after a rescaling of metrics, give the same result.
\end{remark}

\begin{remark} Manifolds with weak holonomy $SU(3)$ structure are (currently) far less plentiful  than their $7$-dimensional counterparts with weak holonomy $G_2$. They are again Einstein manifolds with positive scalar curvature, but the only known examples are $S^6$, the twistor spaces of $S^4$ and $\mathbf{C}P^2$ and $S^3\times S^3$.
\end{remark}

\section{Evolution equations}

We return to a $7$-manifold and consider a fixed cohomology class $\mathcal{A}$ of closed $4$-forms. This is an infinite dimensional affine space, whose tangent space at each point is naturally isomorphic to $\Omega^4_{exact}(M)$. The quadratic form $Q$ defined above provides an indefinite metric on $\mathcal{A}$:
$$(\dot\rho,\dot\rho)=\int_M\gamma\wedge d\gamma$$
where $\dot\rho=d\gamma$.

In this setting the functional  $V(\rho)$ defines a gradient vector field $X$ on $\mathcal{A}$, and the critical points of $V(\rho)$ are the zeros of $X$. More generally, we can consider the gradient flow, and we find:
\begin{theorem}
\label{gradient}
Let $M$ be a closed $7$-manifold and suppose $\rho(t)$ is a closed stable $4$-form which evolves via the gradient flow of the functional $V(\rho)$
restricted to a
cohomology class in $H^4(M,\mathbf{R})$. Then the $4$-form
$$\varphi= dt\wedge \ast \rho + \rho$$
defines a metric with holonomy $Spin(7)$ on the $8$-manifold $N=M\times (a,b)$ for an interval $(a,b)$.
 Conversely, if $N$ is an $8$-manifold with holonomy $Spin(7)$, foliated by equidistant compact hypersurfaces diffeomorphic to $M$, the restriction of the defining $4$-form to each hypersurface evolves as the gradient flow of $V(\rho)$.
  \end{theorem}

\begin{proof} The gradient flow  is  the solution of the equation
$$\frac{dx}{dt}=X$$
where $(X,Y)=D\phi(Y)$ for any vector field $Y$. In our case we have $x$ described by a $4$-form $\rho(t)$ so using the inner product above with $Y=d\gamma$,
$$\left(\frac{dx}{dt},Y\right)=\int_M \frac{\partial \rho}{\partial t}\wedge \gamma=D\phi(Y)=\int_M \ast \rho\wedge d\gamma$$
since $D\phi=\ast \rho$. But this equation holds for all $3$-forms $\gamma$ and thus yields the gradient flow equation
\begin{equation}
\frac{\partial \rho}{\partial t}=d(\ast \rho)
\label{grad}
\end{equation}
But $d\rho=0$, so we obtain
$$d(dt\wedge\ast \rho +\rho)=0$$
and from \cite{Sal} we see that the $4$-form $dt\wedge\ast \rho +\rho$ defines a metric with holonomy $Spin(7)$. In this metric $\Vert dt \Vert^2=1$.

Conversely, if a $Spin(7)$ manifold is foliated by equidistant hypersurfaces, defining the function $t$ to be the distance to  a fixed hypersurface $M$, we can write the defining $4$-form in the form $dt\wedge \ast \rho +\rho$ and the statement that this is closed is equivalent to the gradient flow equation (\ref{grad}).
\end{proof}
\vskip .25cm
We can do something similar with the $6$-dimensional case. Here we take $\mathcal{A}$ to be a cohomology class in $H^3(M,\mathbf{R})$ and $\mathcal{B}$ to be a class in $H^4(M,\mathbf{R})$. Then the tangent space at each point of the product of affine spaces $\mathcal{A}\times\mathcal{B}$ is naturally isomorphic to
$$\Omega^3_{exact}(M)\times\Omega^4_{exact}(M)$$
and (\ref{duality}) shows that the pairing $\langle\rho,\sigma\rangle$ defines formally a symplectic structure on  $\mathcal{A}\times\mathcal{B}$:
$$\omega((\rho_1,\sigma_1),(\rho_2,\sigma_2))=\langle\rho_1,\sigma_2\rangle-\langle\rho_2,\sigma_1\rangle.$$

\begin{remark}
The pairing between  $V$ and $V^*$ to define a symplectic structure on $V\times V^*$ can also be used to define an indefinite metric:
$$((\rho_1,\sigma_1),(\rho_2,\sigma_2))=\langle\rho_1,\sigma_2\rangle+\langle\rho_2,\sigma_1\rangle.$$
The gradient flow of a function of the form $f(\rho)+g(\sigma)$ is then equivalent to the Hamiltonian flow of $f(\rho)-g(\sigma)$. With $V=\Omega^3_{exact}(M)$, the gradient flow shows, as a particular consequence, the relationship between weak holonomy $SU(3)$ and conical $G_2$ metrics as in \cite{Ba}.  Nevertheless, the Hamiltonian  interpretation of the equations has  certain advantages over  the gradient viewpoint, which is why we adopt it in the next theorem.
\end{remark}

\begin{theorem}\label{g2evolve}
Let $\mathcal{A}\in H^3(M,\mathbf{R})$ and $\mathcal{B}\in H^4(M,\mathbf{R})$ be cohomology classes and $(\rho,\sigma)\in \mathcal{A}\times\mathcal{B}$ be stable forms of positive type which evolve via the Hamiltonian flow of the functional $H=V(\rho)-2V(\sigma)$. If for time $t=t_0$, $\rho$ and $\sigma$ satisfy the compatibility conditions $\omega\wedge\rho=0$ and $\phi(\rho)=2\phi(\sigma)$  then the $3$-form
$$\varphi= dt\wedge \omega + \rho$$
(where $\sigma=\omega^2/2$) defines a metric with holonomy $G_2$ on the $7$-manifold $M\times (a,b)$. Conversely, if $N$ is a $7$-manifold with holonomy $G_2$, foliated by equidistant compact hypersurfaces diffeomorphic to $M$, the restriction of the defining closed forms $\rho,\sigma$ to each hypersurface evolves as the Hamiltonian flow of $H$.
\end{theorem}
\begin{proof} As before, we have the derivative of the Hamiltonian $H=V(\rho)-2V(\sigma)$ given by
$$DH(d\alpha,d\beta)=\int_M\hat \rho\wedge d\alpha-\int_M\omega\wedge d\beta$$
and the Hamiltonian vector field $X$ is defined by $\iota(X)\omega=dH$ so in our case we have
$$\int_M\dot\rho\wedge \beta-\int_M\dot\sigma\wedge\alpha=\int_M\hat \rho\wedge d\alpha-\int_M\omega\wedge d\beta$$
and this gives the equations:
\begin{equation}
\frac{\partial\rho}{\partial t}=d\omega
\label{ham1}
\end{equation}
\begin{equation}
\frac{\partial \sigma}{\partial t}=\omega\wedge \frac{\partial\omega}{\partial t}=-d\hat\rho
\label{ham2}
\end{equation}

We shall see that if the compatibility conditions (\ref{compat1}), (\ref{compat2}) between $\omega$ and $\rho$  hold for $t=t_0$, then they hold for all subsequent time.

First consider the condition $\omega\wedge\rho=0$. This  can be viewed as the vanishing of the  moment map for the natural action of $\Diff(M)$ on the symplectic manifold $\mathcal{A}\times \mathcal {B}$.  To see this, for a vector field $X$ on $M$, consider the function
\begin{equation}
\mu_X(\rho,\sigma)=\int_M\iota(X)\sigma\wedge \rho=\int_M\iota(X)\omega\wedge\omega\wedge \rho
\label{moment}
\end{equation}
We have
$$d\mu_X(\dot\rho,\dot\sigma)=\int_M\iota(X)\dot\sigma\wedge\rho+\int_M\iota(X)\sigma\wedge\dot\rho$$
But the first term is
$$-\int_M\dot\sigma\wedge \iota(X)\rho$$
since $\dot\sigma\wedge\rho=0$. Thus
\begin{equation}
d\mu_X(\dot\rho,\dot\sigma)=\int_M\iota(X)\rho\wedge\dot\sigma+\int_M\iota(X)\sigma\wedge\dot\rho
\label{demu}
\end{equation}
Now since $\rho$ and $\sigma$ are closed,
$$\mathcal{L}_X(\rho,\sigma)=(d\iota(X)\rho,d\iota(X)\sigma)$$ so (\ref{demu}) can be written 
$$d\mu_X(\dot\rho,\dot\sigma)=\langle\mathcal{L}_X(\rho,\sigma),(\dot\rho,\dot\sigma)\rangle$$
using the definition of the symplectic form. We deduce that $\mu_X$ is the moment map for $\Diff(M)$, evaluated on $X$. Since $\omega$ is non-degenerate, from (\ref{moment}) $\mu_X$ vanishes for all $X$ if and only if $\omega\wedge\rho=0$. 

Now since the functional $H$ is diffeomorphism invariant, $H$ Poisson commutes with all the functions $\mu_X$. Hence the Hamiltonian flow of $H$ is tangential to the zero set of all the functions $\mu_X$, i.e. the space of pairs $(\rho,\sigma)\in \mathcal{A}\times\mathcal{B}$ such that $\omega\wedge\rho=0$. Thus if $\omega\wedge\rho=0$ holds for $t=t_0$, it holds for all time. Note that it follows then that $\omega\wedge\hat\rho=0$ for all time too.
\vskip .25cm
Next consider the second compatibility condition (\ref{compat2}). The form $\hat\rho$ is defined by the derivative of $\phi$,   $D\phi(\dot\rho)=\hat\rho\wedge\dot\rho$, so the derivative of $\hat\rho$ can be expressed via $D^2\phi$. The volume $\phi(\rho)$ is homogeneous of degree $2$, so its derivative is homogeneous of degree one and hence
$$D^2\phi(\rho,\dot\rho)=D\phi(\dot\rho)=\hat\rho\wedge\dot\rho$$
 We have
\begin{equation}
\frac{\partial \hat\rho }{\partial t}\wedge \rho=D^2\phi(\rho,\frac{\partial \rho }{\partial t})
=D\phi(\frac{\partial \rho }{\partial t})
=\hat\rho \wedge \frac{\partial \rho }{\partial t}
\label{hatderiv}
\end{equation}
So from (\ref{hatderiv})
$$\frac{\partial }{\partial t}(\hat\rho\wedge\rho)=2\hat\rho\wedge\frac{\partial \rho}{\partial t}=2\hat\rho\wedge d\omega$$
from (\ref{ham1}). But $\omega \wedge \hat \rho=0$, and hence $\hat\rho\wedge d\omega=-d\hat\rho\wedge\omega$
Thus
$$\frac{\partial }{\partial t}(\hat\rho\wedge\rho)=-2d\hat\rho\wedge\omega=2\omega^2\frac{\partial \omega}{\partial t}=\frac{2}{3}\frac{\partial}{\partial t}\omega^3$$
and so if $\phi(\rho)=2\phi(\sigma)$ at $t=t_0$ then it holds for all $t$.
\vskip .25cm
The evolution equation thus preserves the $SU(3)$ geometry on $M$, and we define the $3$-form
$$\varphi=dt\wedge\omega+\rho$$
From (\ref{ham1}), $d\varphi=0$. From (\ref{ast}), using the $G_2$ metric defined by $\varphi$, we have
$$\ast\varphi=dt\wedge\hat\rho-\sigma$$
and from (\ref{ham2}), $d\ast\varphi=0$. Thus $\varphi$ defines a metric of holonomy $G_2$ on $M\times (a,b)$ for some interval.
\end{proof}

\section{Examples}

It is a well-known principle that if we are looking for a $G$-invariant critical point of an invariant function, we need only consider critical points of the same function restricted to the fixed point set of $G$. Thus if $G$ acts on $M$, a $G$-invariant special metric of any of the types we are considering can be found by restricting the volume functionals to $G$-invariant forms. This applies in particular when $M$ is homogeneous under $G$, in which case we are reduced to a finite-dimensional variational problem.

Since metrics of holonomy $G_2$ and $Spin(7)$ have zero Ricci tensor, we have no interesting compact homogeneous examples. We can nevertheless apply the principle to compact spaces for  weak holonomy $G_2$ and $SU(3)$, and also to the evolution equation to find non-compact examples with holonomy $Spin(7)$ and $G_2$.

\subsection{A $7$-dimensional example}

We consider $M^7=S^7$ as a principal $SU(2)$ bundle over $S^4$, and $G=SO(5)$ acting transitively on it with stabilizer $SU(2)$. All the relevant cohomology classes are trivial here, so we are dealing with exact forms, which must be built out of  the three invariant $1$-forms $\alpha_1,\alpha_2,\alpha_3$ which are the components of the $SU(2)$-connection form, and the components $\omega_1,\omega_2,\omega_3$ of the
curvature. These are defined by:
\begin{equation}
d\alpha_1+2\alpha_2\alpha_3=\omega_1
\label{curvature}
\end{equation}
and similar expressions for $\alpha_2,\alpha_3$. The curvature forms themselves satisfy the Bianchi identity
\begin{equation}
d\omega_1=2(\omega_2\alpha_3-\omega_3\alpha_2)\,\,etc.
\label{bianchi}
\end{equation}
A basis for the exact invariant $4$-forms is provided by
$$u_1=d(\alpha_1\omega_1), u_2=d(\alpha_2\omega_2),u_3=d(\alpha_3\omega_3),u_4=d(\alpha_1\alpha_2\alpha_3)$$
and the quadratic form
$$Q(d\gamma)=(d\gamma,d\gamma)=\int_M\gamma d\gamma$$
defines a metric on this space. We have, for example,
$$(u_1,u_1)=\int_M\alpha_1\omega_1d(\alpha_1\omega_1)=\int_M\alpha_1\omega_1(-2\alpha_2\alpha_3)\omega_1=-2\int_M\alpha_1\alpha_2\alpha_3\omega_1^2$$
using (\ref{curvature}) and (\ref{bianchi}). Here, because the connection is anti-self-dual, we have
$$\omega_1^2=\omega_2^2=\omega_3^2=\nu,\quad \omega_1\omega_2=\omega_2\omega_3=\omega_3\omega_1=0.$$
Continuing this way, and normalizing the integral of $\alpha_1\alpha_2\alpha_3\nu$, we find the metric with respect to this basis to be
$$(u_1,u_1)=(u_2,u_2)=(u_3,u_3)=2,\quad (u_1,u_2)=(u_2,u_3)=(u_3,u_1)=-2$$
$$(u_1,u_4)=(u_2,u_4)=(u_3,u_4)=1,\quad (u_4,u_4)=0.$$

We now need the volume $V(\rho)$ for a general exact invariant  $4$-form
$$\rho=\sum_1^3x_id(\alpha_i\omega_i)+2x_4d(\alpha_1\alpha_2\alpha_3)$$
which from (\ref{curvature}) and (\ref{bianchi}) can be written 
$$\rho=(x_1+x_2+x_3)\nu+k_1\alpha_2\alpha_3\omega_1+k_2\alpha_3\alpha_1\omega_2+k_3\alpha_1\alpha_2\omega_3$$
where $k_1=2(x_4-x_1+x_2+x_3)$ etc.
Instead of working out the volume as in (8.4), it is  easier to transform $\rho$ to the standard form (\ref{ast}). In fact we can rewrite that formula as
$$\ast\varphi=e_1e_2e_3e_4+e_5e_6(e_2e_3+e_4e_1)+e_6e_7(e_4e_3+e_1e_2)+e_7e_5(e_1e_3+e_2e_4)$$
Then putting $e_5=y_1\alpha_1, e_6=y_2\alpha_2, e_7=y_3\alpha_3$ and rescaling $e_i=y_4v_i$ for $1\le i\le 4$, where $v_i$ is a local orthonormal basis such that $\omega_1=(v_4v_3+v_1v_2)$ etc. we can transform $\rho$ into $\ast\varphi$ if 
$$k_1=y_2y_3y_4^2,\quad k_2=y_3y_1y_4^2,\quad k_3=y_1y_2y_4^2,\quad y_4^4=-2(x_1+x_2+x_3)$$
 The Riemannian volume is then
\begin{equation}
V(\rho)=y_1y_2y_3y_4^4
\label{voleq}
\end{equation}
Using the expression for the inner products above, we find the gradient flow equation
$$2(-\dot x_1+\dot x_2+\dot x_3+\dot x_4)=\frac{\partial V}{\partial x_1},\,\, etc.$$
$$2(\dot x_1+\dot x_2+\dot x_3)=\frac{\partial V}{\partial x_4}$$
Converting to coordinates $y_i$ defined by
$$y_1y_2y_4^2=2(x_1+x_2-x_3+x_4)\,\,etc., \quad y_4^4=-2(x_1+x_2+x_3),$$
we find the equations
\begin{eqnarray*}
\dot y_1&=&-1+\frac{1}{2y_2y_3}(y_2^2+y_3^2-y_1^2)+\frac{y_1^2}{2y_4^2}\\
\dot y_2&=&-1+\frac{1}{2y_3y_1}(y_3^2+y_1^2-y_2^2)+\frac{y_2^2}{2y_4^2}\\
\dot y_3&=&-1+\frac{1}{2y_1y_2}(y_1^2+y_2^2-y_3^2)+\frac{y_3^2}{2y_4^2}\\
4y_4\dot y_4&=&-(y_1+y_2+y_3)
\end{eqnarray*}
Putting $-\alpha=y_2=y_3,\beta=-y_1,\gamma=y_4$ gives the equations solved by \v Cvetic et al in \cite{C}. Their solution represents a $Spin(7)$ manifold with  an extra $S^1$  symmetry.
\vskip .25cm
The constrained variational problem to get a weak holonomy $G_2$ structure consists of finding critical points of $V$ subject to the condition that the quadratic form is constant. This is precisely when the gradient field is parallel to the position vector. As an example, if $y_1=y_2=y_3=y$ in the equations above, then we obtain
$$\lambda y=-\frac{1}{2}+\frac{y^2}{2y_4^2},
\quad 4\lambda y_4^2=-3y$$
which yields $y=-3/10\lambda$, and this is the squashed $7$-sphere.
\vskip .25cm
In this same symmetric situation the gradient flow equations become
$$\dot y=-\frac{1}{2}+\frac{y^2}{y_4^2},\quad 2y_4\dot y_4=-\frac{3}{2}y$$
Putting $s=y_4^2$ the equation becomes
$$y\frac{dy}{ds}+\frac{y^2}{3s}=\frac{1}{3}$$
which is solved easily:
$$y^2=\frac{2}{5}s+cs^{-2/3}.$$
When $c=0$ we have the cone on the squashed $S^7$ and for $c\ne 0$, the complete Bryant-Salamon metric on the spin bundle over $S^4$ \cite{Sal}.

\subsection{A $6$-dimensional example}

We take now $M^6=S^3\times S^3$ and the forms to be invariant under $G$ given by the left action of $S^3\times S^3$ and \textit{anti-invariant} under the $\mathbf{Z}/2$ action which interchanges the two factors.
We have left-invariant $1$-forms $$\sigma_1,\sigma_2,\sigma_3,\quad \Sigma_1,\Sigma_2,\Sigma_3$$ on the two factors, satisfying
$$d\sigma_1=-\sigma_2\sigma_3,\,\,etc.$$

We shall take  $\mathcal{A}\in H^3(S^3\times S^3)$ to be $\nu_1-\nu_2$ where $\nu_1$ is the pull-back of a  generator of $H^3(S^3)$ on the first factor and $\nu_2$ the pull-back from the second. This is clearly anti-invariant under the $\mathbf{Z}/2$ action. Since $H^4(S^3\times S^3)=0$, $\mathcal{B}$ must be the trivial class.

We find the following two three-dimensional spaces of closed invariant $3$-forms $\rho$ and $4$-forms $\sigma$:
\begin{equation}
\rho=\sigma_1\sigma_2\sigma_3- \Sigma_1\Sigma_2\Sigma_3+x_1(\sigma_1\Sigma_2\Sigma_3-\sigma_2\sigma_3\Sigma_1)+\cdots
\label{rhoform}
\end{equation}
\begin{equation}
\sigma=y_1\sigma_2\Sigma_2\sigma_3\Sigma_3+y_2\sigma_3\Sigma_3\sigma_1\Sigma_1+y_3\sigma_1\Sigma_1\sigma_2\Sigma_2
\label{sigmaform}
\end{equation}
Since $d(\sigma_1\Sigma_1)=\sigma_1\Sigma_2\Sigma_3-\sigma_2\sigma_3\Sigma_1$ we see that $\rho$ lies in the fixed cohomology class $\mathcal{A}$. The $4$-form $\sigma$ is clearly closed. Moreover $\omega\wedge\rho=0$ for all of these forms. We can work out the pairing between invariant exact $3$-forms and invariant exact $4$-forms in a straightforward manner, e.g.
$$\langle \sigma_2\Sigma_2\sigma_3\Sigma_3, d(\sigma_1\Sigma_1)\rangle=\int_{S^3\times S^3}\sigma_2\Sigma_2\sigma_3\Sigma_3\sigma_1\Sigma_1=\nu_1[S^3]\nu_2[S^3]\ne 0.$$
It turns out that the symplectic form is a multiple of
$$dx_1\wedge dy_1+ dx_2\wedge dy_2+dx_3\wedge dy_3$$
We can use (8.2) to find $V(\rho)$ and then we find
$$V(\rho)^2=(1+x_1+x_2+x_3)(x_2+x_3-x_1-1)(x_3+x_1-x_2-1)(x_1+x_2-x_3-1)$$
$$V(\sigma)^2=y_1y_2y_3$$
Since we are only interested in solutions where $V(\rho)=2V(\sigma)$, after a change of parameter $t$, it is equivalent to consider the Hamiltonian flow of $4V^2(\sigma)-V^2(\rho)$ i.e. for the Hamiltonian
$$H=4y_1y_2y_3-(1+x_1+x_2+x_3)(x_2+x_3-x_1-1)(x_3+x_1-x_2-1)(x_1+x_2-x_3-1)$$
The six first order equations arising from this are equivalent to those in \cite{Br}, where the authors produce an explicit complete solution with $x_2=x_3$ and $y_2=y_3$. They express their solution using an orthonormal basis $A_j(\sigma_j-\Sigma_j), B_j(\sigma_j+\Sigma_j)$ for the metric and our coordinates relate to theirs by $x_1=A_1A_2A_3+A_1A_2B_3+A_3B_1B_2-A_1B_2B_3$ etc. and $y_1=4A_2B_2A_3B_3$ etc.
\vskip .25cm
When $x_1=x_2=x_3=x$ and $y_1=y_2=y_3=y$, we don't need to solve the equations since the solutions are given by the vanishing of the Hamiltonian:
$$4y^3-(1+3x)(x-1)^3=0$$
This is the metric originally found by Bryant and Salamon \cite{Sal} on the spin bundle over $S^3$.
\vskip .25cm
For a weak holonomy $SU(3)$ metric we need the   cohomology classes $\mathcal{A},\mathcal{B}$ to be trivial, and then  the ``1" terms disappear in the functional. In the fully symmetric situation, we need then to  find critical points of $8y^{3/2}+3\sqrt{3}x^2$  subject to the condition $xy=c$. The solution $y^{7/2}=\sqrt{3}c^2/2$ gives the weak holonomy $SU(3)$ metric on $S^3\times S^3$.

\section{Appendix: Definition of the volumes}
We assume $V$ to be oriented and then make all volume forms positive.
\subsection{$\mathbf{p=2m-2,n=2m}$}

The vector space $V$ is $2m$-dimensional and $\rho\in
\Lambda^{2m-2}V^*$. Use the isomorphism $$\Lambda^{2m-2}V^*\cong
\Lambda^2V\otimes\Lambda^{2m}V^*$$ to write $\rho$ as $\sigma\in
\Lambda^2V\otimes\Lambda^{2m}V^*$. Then we have using the  exterior
product $$\sigma^{2m}\in
\Lambda^{2m}V\otimes(\Lambda^{2m}V^*)^m\cong
(\Lambda^{2m}V^*)^{m-1}$$ We define
$$\phi(\rho)=\vert\sigma^{2m}\vert^{1/(m-1)}.$$

\subsection {$\mathbf{n=6,p=3}$}

We have a $6$-dimensional space $V$ and $\rho\in \Lambda^3V^*$.
For $v\in V$ define $$K(v)=\iota(v)\rho\wedge\rho\in
\Lambda^5V^*\cong V\otimes \Lambda^6V^*.$$ Then $$\tr(K^2)\in
(\Lambda^6V^*)^2.$$ The stable forms with stabilizer
$SL(3,\mathbf{C})$ are characterized by $\tr(K)^2<0$ and there
$$\phi(\rho)=\vert\sqrt{-\tr K^2}\vert\in \Lambda^6V^*.$$

\subsection{$\mathbf{n=7,p=3}$}

The space $V$ is $7$-dimensional and $\rho\in\Lambda^3V^*$. Given
$v,w\in V$ form $$\iota(v)\rho\wedge\iota(w)\rho\wedge\rho\in
\Lambda^7V^*.$$ This is a symmetric bilinear form on $V$ with
values in $\Lambda^7V^*$, and so a linear map $$G:V\rightarrow
V^*\otimes \Lambda^7V^*$$ We have $$\det G\in (\Lambda^7V^*)^9$$
and we define $$\phi(\rho)=\vert\det G\vert^{1/9}.$$

\subsection{$\mathbf{n=7,p=4}$}

Take $\rho\in \Lambda^4V^*$ but identify with an element
$\sigma\in \Lambda^3V\otimes\Lambda^7V^*$. Now apply the
construction of the previous case and we find a linear map
$$H:V^*\rightarrow V\otimes (\Lambda^7V^*)^2$$ This gives $$\det
H\in (\Lambda^7V^*)^{12}$$ and we define $$\phi(\rho)=\vert\det
H\vert^{1/12}.$$

\subsection{$\mathbf{n=8,p=3}$ and $\mathbf{p=5}$}

The space $V$ is $8$-dimensional and $\rho\in \Lambda^3V^*$. For
$v,w\in V$ form $$\iota(v)\rho\wedge\iota(w)\rho\wedge\rho\in
\Lambda^7V^*\cong V\otimes\Lambda^8V^*$$ This gives an element
$$d\in V^*\otimes V^*\otimes V\otimes \Lambda^8V^*$$ Take
$d\otimes d$ and contract to get $$G\in V^*\otimes V^*\otimes
(\Lambda^8V^*)^2$$ or equivalently $$G:V\rightarrow V^*\otimes
(\Lambda^8V^*)^2.$$Then we have $$\det G \in (\Lambda^8V^*)^{18}$$
and define $$\phi(\rho)=\vert\det G\vert^{1/18}.$$ For $p=5$ we use
$\Lambda^5V^*\cong \Lambda^3V\otimes \Lambda^8V^*$ and proceed
similarly.

\bibliographystyle{amsalpha}

\end{document}